\documentclass[12pt]{amsart}

\usepackage[latin1]{inputenc}
\usepackage[T1]{fontenc}
\usepackage{amsmath,amssymb,amsthm}
\usepackage{graphicx}
\usepackage{epsfig}
\usepackage{hyperref}

\usepackage[square,authoryear,sort&compress]{natbib}
\usepackage[total={8.5in,11in},text={6.25in,8.25in}]{geometry}

\begin{document}

\newtheorem{theorem}{Theorem}[section]
\newtheorem{lemma}[theorem]{Lemma}
\newtheorem{corollary}[theorem]{Corollary}
\newtheorem{proposition}[theorem]{Proposition}
\theoremstyle{definition}
\newtheorem{definition}[theorem]{Definition}
\newtheorem{example}[theorem]{Example}
\newtheorem{claim}[theorem]{Claim}
\newtheorem{xca}[theorem]{Exercise}

\theoremstyle{remark}
\newtheorem{remark}[theorem]{Remark}

%    Blank box placeholder for figures (to avoid requiring any
%    particular graphics capabilities for printing this document).
\newcommand{\blankbox}[2]{%
  \parbox{\columnwidth}{\centering
%    Set fboxsep to 0 so that the actual size of the box will match the
%    given measurements more closely.
    \setlength{\fboxsep}{0pt}%
    \fbox{\raisebox{0pt}[#2]{\hspace{#1}}}%
  }%
}

\bibliographystyle{plain}

\title{The neighborhood complex of a random graph}

\author{Matthew Kahle}
\address{Department of Mathematics, University of Washington}
\email{kahle@math.washington.edu}

\date{January 31,2006}

\thanks{Supported by NSF grant DMS-9983797.}

\maketitle

%\vspace*{-0.5in}

\begin{abstract}
\noindent For a graph $G$, the neighborhood complex
$\mathcal{N}[G]$ is the simplicial complex having all subsets of
vertices with a common neighbor as its faces. It is a well known
result of Lov\'{a}sz that if $\| \mathcal{N}[G] \|$ is
$k$-connected, then the chromatic number of $G$ is at least $k +
3$.

We prove that the connectivity of the neighborhood complex of a
random graph is tightly concentrated, almost always between $1/2$
and $2/3$ of the expected clique number. We also show that the
number of dimensions of nontrivial homology is almost always
small, $O(\log{d})$, compared to the expected dimension $d$ of the
complex itself.

\end{abstract}

\section{Introduction}

In 1978, L\'{a}szl\'{o} Lov\'{a}sz proved Kneser's conjecture
\cite{Lovasz}, that if the $n$-subsets of a $(2n+k)$-set are
partitioned into $k+1$ families, at least one family contains a
disjoint pair. He restated the problem graph theoretically and
then proved a more general theorem about graph coloring.

All our graphs will be simple undirected graphs, with no loops or
multiple edges. For a graph $G$, a {\it $k$-coloring} is a
function $f:V(G) \rightarrow \{ 1, 2, \ldots, k \}$ such that
$f(x) \neq f(y)$ whenever $\{ x, y\} \in E(G)$. The {\it chromatic
number} $\chi(G)$ is the minimum $k$ such that $G$ admits a
$k$-coloring. Define the {Kneser graph} $KG(n,k)$ to have all
$n$-subsets of a $(2n+k)$-set as its vertices, with edges between
disjoint pairs. The Kneser conjecture is equivalent to the claim
that $\chi(KG(n,k)) \ge k+2$.

There are several simplicial complexes naturally assoicated with a
graph $G$. The clique complex $X(G)$ is the simplicial complex
$X(G)$ on vertex set $V(G)$ whose simplices are all complete
subgraphs of $G$. In another article \cite{Kahle}, we study the
clique complex of a random graph. The {\it neighborhood complex}
$\mathcal{N}[G]$ is the simplicial complex on $V(G)$ which has all
subsets of $V(G)$ with a common neighbor for its faces. For
example, the neighborhood complex of the complete graph $K_n$ has
all proper subsets of the vertices for its faces, so it is the
boundary of an $(n-1)$-dimensional simplex. Its geometric
realization $\| \mathcal{N}[K_n] \|$ is homeomorphic to an
$(n-2)$-dimensional sphere $\mathbb{S}^{n-2}$. (We denote the
geometric realization of a simplicial complex $\Delta$ by $\|
\Delta \|$.)

A topological space $X$ is said to be {\it $k$-connected} if every
map from a sphere $\mathbb{S}^{n} \rightarrow X$ extends to a map
from the ball $\mathbb{B}^{n+1} \rightarrow X$ for $n=0, 1,
\ldots, k$. On the way to proving the Kneser conjecture,
Lov\'{a}sz proved the following \cite{Lovasz}.

\begin{theorem}[Lov\'{a}sz] \label{Theorem A} If $\| \mathcal{N}[G] \|$ is $k$-connected, then $\chi(G) \ge k+3$.
\end{theorem}

In the case of the Kneser graphs this lower bound is tight,
matching an easy upper bound to give an exact answer. It seems
natural to ask how good is the bound for ``typical'' graphs. For
our purposes, a typical graph is the random graph $G(n,p)$
\cite{Bollo}.

The random graph $G(n,p)$ is the probability space of all graphs
on a vertex set of size $n$ with each edge inserted independently
with probability $p$. Frequently, one considers $p$ to be a
function of $n$ and asks whether the graph is likely to have some
property as $n \rightarrow \infty$. We say that $G(n,p)$ {\it
almost always (a.a.)} has property $\mathcal{P}$ if
$\mbox{Prob}[G(n,p) \in \mathcal{P}] \rightarrow 1$ as $n
\rightarrow \infty$. The main goal of this article is to
understand some of the most basic topological features of $\|
\mathcal{N}[G(n,p)] \|$.

\section{Statement of results}

Let $p=p(n)$ be a monotone function of $n$, and let $i,j,k,$ and
$l$ be integer valued monotone functions of $n$. In the asymptotic
notation that follows, $n \rightarrow \infty$ is the free
variable. Homology is understood to be reduced with coefficients
in $\mathbb{Z}$ throughout.

\begin{theorem}\label{Theorem 1} If ${n \choose i} (1-p^i)^{n-i} = o(1)$ then
 $\mathcal{N}[G(n,p)]$ is a.a. $(i-2)$-connected.
\end{theorem}

\begin{theorem}\label{Theorem 2} If ${n \choose j}{n \choose k} p^{jk} = o(1)$ then
$\mathcal{N}[G(n,p)]$ a.a. strong deformation retracts to a
simplicial complex of dimension at most $j+k-3$.
\end{theorem}

In proving Theorem \ref{Theorem 2}, we make use of the following
lemma, which might be of independent interest.

\begin{lemma}\label{lemma 1} If $H$ is any graph not containing a
complete bipartite subgraph $K_{a,b}$ then $|| \mathcal{N}[H] ||$
strong deformation retracts to a complex of dimension at most
$a+b-3$.
\end{lemma}

A $d$-connected complex has trivial homology through dimension $d$
by the Hurewicz theorem. Also, $\widetilde{H}_k(\Delta)=0$
whenever $k$ is greater than the dimension of $\Delta$, and strong
deformation retracts are homotopy equivalences and preserve
homology. Hence Theorems \ref{Theorem 1} and \ref{Theorem 2} bound
the possible dimensions of nontrivial homology from below and
above. We give special cases of the theorems as corollaries. First
fix $p=1/2$, where $G(n,p)$ is the uniform distribution on all
graphs on vertex set $[n]$.

\begin{corollary} \label{corollary 1} If $p=1/2$ and $\epsilon > 0$ then a.a. $\widetilde{H}_l(\|
\mathcal{N}[G(n,p)] \|)=0$ for $l \le (1-\epsilon) \log_2{n}$ and
$l \ge (4+\epsilon) \log_2{n}$.
\end{corollary}

Note for comparison that the dimension of the neighborhood complex
is one less than the maximum vertex degree, so when $p=1/2$ we
expect it to be slightly more than $n/2$. Next, fix $l$ and check
that $\widetilde{H}_l$ is trivial outside a certain range of $p$.

\begin{corollary} \label{corollary 2} Let $p=n^\alpha$ with $\alpha \in [-2,0]$.
If $\alpha  > \frac{-1}{l+2}$ then a.a. $\widetilde{H}_l(\|
\mathcal{N}[G(n,p)] \|)=0$. For $l$ even, if $\alpha <
\frac{-4}{l+2}$ then a.a. $\widetilde{H}_l(\| \mathcal{N}[G(n,p)]
\|)=0$. For $l$ odd, if $\alpha < \frac{-4(l+2)}{(l+1)(l+3)}$ then
a.a. $\widetilde{H}_l(\| \mathcal{N}[G(n,p)] \|)=0$.
\end{corollary}

For a partial converse to Corollaries \ref{corollary 1} and
\ref{corollary 2}, we exhibit explicit nontrivial homology classes
by retracting onto random spheres. Recall that a {\it clique} of
order $n$ is a complete subgraph on $n$ vertices.

\begin{definition}
Let the graph $X_n$ have vertex set $\{ u_1, u_2, \ldots, u_n \}
\coprod \{ v_1, v_2, \ldots, v_n \}$, such that $\{ u_1, u_2,
\ldots, u_n \}$ spans a clique, and $u_i$ is adjacent to $v_j$
whenever $i \neq j$.
\end{definition}

\begin{theorem}\label{Theorem 3} If $H$ is any graph containing a maximal clique of
order $n$ that can't be extended to an $X_n$ subgraph, then $\|
\mathcal{N}[H] \|$ retracts onto a sphere $\mathbb{S}^{n-2}$.
\end{theorem}

\begin{corollary} \label{corollary 3} If $p=1/2$, $\epsilon > 0$,
and $(4/3+\epsilon) \log_2{n} <k < (2-\epsilon) \log_2{n}$, then
a.a. $\widetilde{H}_k(\| \mathcal{N}[G(n,p)] \|) \neq 0$.
\end{corollary}

\begin{corollary} \label{corollary 4}
Let $p=n^\alpha$ with $\frac{-2}{k+1}< \alpha <
\frac{-4}{3(k+1)}$, then a.a. $ \widetilde{H}_k(\|
\mathcal{N}[G(n,p)] \|) \neq 0$.
\end{corollary}

\section{Proofs}

We first prove that if ${n \choose i} (1-p^i)^{n-i} = o(1)$ then
 $\mathcal{N}[G(n,p)]$ is a.a. $(i-2)$-connected.

 \begin{proof}[Proof of Theorem \ref{Theorem 1}] A simplex complex is {\it $i$-neighborly} if
 every $i$ vertices span a face. By simplicial approximation, if
 a complex is $i$-neighborly then it is $(i-2)$-connected. The
 probability that a given set of $i$ vertices in $G(n,p)$ has no common
 neighbor is $(1-p^i)^{n-i}$. Then the total probability that any
 set of $i$ vertices doesn't have a common neighbor is bounded above by ${n \choose i}
 (1-p^i)^{n-i}=o(1)$. So a.a. every such set has some common neighbor,
 hence spans a face in the neighborhood complex. So $\mathcal{N}[G(n,p)]$
 is a.a. $i$-neighborly and $(i-2)$-connected.
 \end{proof}

Next we prove that if $H$ is any graph not containing a complete
bipartite subgraph $K_{a,b}$ then $|| \mathcal{N}[H] ||$ strong
deformation retracts to a complex of dimension at most $a+b-3$.

\begin{proof}[Proof of Lemma \ref{lemma 1}]

For a poset $Q$, the {\it order complex} $\Delta(Q)$ is the
simplicial complex  of all chains in $Q$. For a simplicial complex
$S$, let $P(S)$ denote its face poset. To avoid proliferation of
notation, we denote the geometric realization of the order complex
of a poset $Q$ by $\|Q\|$ rather than $\|\Delta(Q)\|$.

For a vertex $x$ of $H$, let $\Gamma(x)$ denote the set of common
neighbors of $x$. Similarly, for any face in the neighborhood
complex  $X \in \mathcal{N}[H]$, let $\Gamma(X)=\cap_{x \in X}
\Gamma(x) $. (This map is used in Lov\'{a}sz's paper
\cite{Lovasz}.) Note that $\Gamma$ is an order reversing self-map
of $P(\mathcal{N}[H])$, abbreviated for the rest of this proof by
$P$. So we can define an order preserving poset map $v:P
\rightarrow P$ by $v(X)=\Gamma(\Gamma(X))$. It's also easy to
check that $\Gamma^3=\Gamma$, so $\Gamma^4=\Gamma^2$, and $v^2=v$.
Since $v(X) \supseteq X$ for every $X$, a standard theorem in
combinatorial homotopy theory \cite{Bjorner} gives that $\| v(P)
\|$ is a strong deformation retract of $\|P  \|$.

An $m+1$-dimensional face in $v(P)$ is a chain of faces in
$\mathcal{N}[H]$, $X_1 \subsetneq X_2 \subsetneq \cdots \subsetneq
X_{m+2}$. Set $Y_i = \Gamma (X_i)$ and we have $Y_1 \supsetneq Y_2
\supsetneq \cdots \supsetneq Y_{m+2}$. (The inclusions $Y_i
\subsetneq Y_{i+1}$ are strict since $\Gamma^3=\Gamma$.)

Suppose $a+b-3 \le m$. Since the inclusions $X_i \subsetneq
X_{i+1}$ are strict and $X_1$ is nonempty, $X_a$ contains at least
$a$ vertices. Similarly, $Y_a$ contains at least $m+3-a \ge b$
vertices. But $Y_p = \Gamma(X_p)$, so $X_a$ and $Y_b$ span a
complete bipartite subgraph $K_{a,b}$. Then if the dimension of
$v(P)$ is at least $m+1$, $H$ contains $K_{a,b}$ subgraphs for
every $a$ and $b$ such that $a+b-3 \le m$, which is the claim.

\end{proof}

Now we apply Lemma \ref{lemma 1} to check that if ${n \choose j}{n
\choose k} p^{jk} = o(1)$ then $\mathcal{N}[G(n,p)]$ a.a. strong
deformation retracts to a simplicial complex of dimension at most
$j+k-3$.

\begin{proof}[Proof of Theorem \ref{Theorem 2}]

Let $U$ and $V$ be vertex subsets of $G(n,p)$ of order $j$ and $k$
respectively. The probability that they span a complete bipartite
graph with parts $U$ and $V$ is $p^{jk}$. So the total probability
that there are any $K_{j,k}$ subgraphs is bounded above by ${n
\choose j}{n \choose k} p^{jk} = o(1)$. There are a.a. no such
subgraphs, so the claim follows by Lemma \ref{lemma 1}.

\end{proof}

Checking Corollaries \ref{corollary 1} and \ref{corollary 2} is
now a straightforward computation.

\begin{proof}[Proof of Corollary \ref{corollary 1}] Let $p=1/2$ and $\epsilon > 0$. If $l \le (1-\epsilon) \log_2{n}$,
then ${n \choose l} (1-p^l)^{n-l}$

$$={n \choose l} (1-(1/2)^{l})^{n-l}$$

$$\le n^l e^{-(1/2)^l (n-l)}$$

$$\le n^{\log_2{n}} e^{-n^{-1+\epsilon}
(n-\log_2{n})} $$

$$= \exp{(\log{n}\log_2{n} - n^{\epsilon}+n^{-1+\epsilon}\log_2{n})}$$

$$ = o(1).$$

Then Theorem \ref{Theorem 1} gives that $\mathcal{N}[G(n,p)]$ is
a.a. $(l-2)$-connected. Since $\epsilon$ doesn't appear anywhere
in the conclusion of the theorem, we can replace it be a slightly
smaller $\epsilon$ and for large enough $n$, $l+2 \le (1-\epsilon)
\log_2{n}$ and this gives that $\mathcal{N}[G(n,p)]$ is a.a.
$l$-connected.

On the other hand, suppose $l \ge (4+\epsilon) \log_2{n}$ and let
$j = \lceil l/2 \rceil$.

$${n \choose j}^2 \left( \frac{1}{2} \right) ^{j^2}$$

$$\le n^{2(2+\epsilon/2)\log_2{n}} \left( \frac{1}{2} \right)^{(2+\epsilon/2)^2 (\log_2{n})^2 }$$

$$= n^{2(2+\epsilon/2)\log_2{n}-(2+\epsilon/2)^2 \log_2{n}}$$

$$= o(n^{-\epsilon\log_2{n}})$$

$$=o(1).$$

Then Theorem \ref{Theorem 2} gives that $\mathcal{N}[G(n,p)]$ a.a.
strong deformation retracts to a complex of dimension at most
$2j-3 \le l-1$.

\end{proof}

\begin{proof}[Proof of Corollary \ref{corollary 2}] Let $p=n^\alpha$ and suppose first that
$\alpha  > \frac{-1}{l+2}$.

$${n \choose l+2} (1-p^{l+2})^{n-l}$$

$$\le n^{l+2} e^{-n^{\alpha (l+2)} (n-l)}$$

$$\le \exp{[(l+2) \log{n} -n^{1+\alpha(l+2)} +  n^{\alpha(l+2)}l] }$$

$$=o(1),$$

\noindent since $l$ is constant and $1+\alpha(l+2)>0$. Then
Theorem \ref{Theorem 1} gives that $\mathcal{N}[G(n,p)]$ is a.a.
$l$-connected. So a.a. $\widetilde{H}_l(\| \mathcal{N}[G(n,p)]
\|)=0$.

Now suppose $l$ is even and $\alpha < \frac{-4}{l+2}$. Set $j =
\frac{l+2}{2}$.

$${n \choose j}^2 (n^\alpha) ^{j^2}$$

$$\le n^{l+2} n^{\alpha(l+2)^2 /4}$$

$$= n^{(l+2)(4+\alpha(l+2))/4}$$

$$=o(1),$$

\noindent since $4+\alpha(l+2) < 0$. So Theorem \ref{Theorem 2}
gives that $\mathcal{N}[G(n,p)]$ a.a. strong deformation retracts
to a complex of dimension at most $2j-3=l-1$.

Similarly, suppose $l$ is odd and $\alpha <
\frac{-4(l+2)}{(l+1)(l+3)}$. Set $j = \frac{l+1}{2}$.

$${n \choose j}{n \choose j+1} (n^\alpha)^{ j(j+1)}$$

$$\le n^{2j+1} n^{\alpha j(j+1)}$$

$$\le n^{l+2} n^{\alpha \frac{l+1}{2} \frac{l+3}{2} }$$

$$=o(1).$$

Then Theorem \ref{Theorem 2} gives that $\mathcal{N}[G(n,p)]$ a.a.
strong deformation retracts to a complex of dimension at most
$2j-2 \le l-1$. In both the even and odd cases $\widetilde{H}_l(\|
\mathcal{N}[G(n,p)] \|)=0$.

\end{proof}

Recall that the graph $X_n$ has vertex set $\{ u_1, u_2, \ldots,
u_n \} \coprod \{ v_1, v_2, \ldots, v_n \}$, such that $\{ u_1,
u_2, \ldots, u_n \}$ spans a clique, $\{ v_1, v_2, \ldots, v_n \}$
is an independent set, and $u_i$ is adjacent to $v_j$ whenever $i
\neq j$. We show now that if $H$ is any graph containing a maximal
clique $\{ u_1, u_2, \ldots, u_n \}$ that isn't contained in an
$X_n$ subgraph, then $\| \mathcal{N}[H] \|$ retracts onto a sphere
$\mathbb{S}^{n-2}$.

\begin{proof} [Proof of Theorem \ref{Theorem 3}] Suppose $H$ contains
 a clique $X=\{ u_1, u_2, \ldots, u_n
\}$ that isn't contained in any larger clique or $X_n$ subgraph.
The induced subcomplex of $\mathcal{N}[G(n,p)]$ on $X$ is a
topological sphere $\mathbb{S}^{n-2}$, since $X$ itself is not a
face by assumption of maximality of the clique. Define a map on
vertices $r_1: V(H) \rightarrow V(H)$ by $r_1(x)=x$ for $x \in X$
and $r_1(x)=u_1$ otherwise.

The only possible obstruction to $r_1$ extending to a simplicial
map $\widetilde{r_1}: \mathcal{N}[H] \rightarrow \mathcal{N}[H]$,
is an $(n-1)$-dimensional face getting mapped onto $X$. This
happens only if for some vertex $u_1^*$, $X \cup \{ u_1^* \} -
u_1$ has a common neighbor $v_1$. Note that $u_1$ isn't adjacent
to $v_1$ since then $X \cup \{ v_1 \}$ would be an extension of
$X$ to a larger clique. Similarly, replacing $u_1$ with $u_i$ for
$i=2, \ldots, n$. If none of the candidate maps $r_1, \ldots, r_n$
extends to a simplicial map, then the $v_i$ are clearly distinct,
since the $u_i$ are distinct and $u_i$ is adjacent to $v_j$ if and
only if $i \neq j$. But this yields an $X_n$ subgraph containing
$X$. Otherwise $\| \mathcal{N}[H] \|$ retracts onto $\| X \| =
\mathbb{S}^{n-2}$ as claimed, via one of these maps.

\end{proof}

\begin{proof} [Proof of Corollary \ref{corollary 3}] Let $p=1/2$ and
$\epsilon>0$. It is well known that $G(n,p)$ a.a. contains maximal
cliques of every order $k$ with $(1+\epsilon) \log_2{n} <k <
(2-\epsilon) \log_2{n}$ \cite{Bollo}. We need only check that
there are a.a. no $X_k$ subgraphs when $k>(4/3+\epsilon)
\log_2{n}$. Note that $X_k$ has $2k$ vertices and $3k(k-1)/2$
edges. Then the probability that $G(n,p)$ contains a copy of $X_k$
is bounded above by

$$(2k)! {n \choose 2k} \left( \frac{1}{2} \right) ^{3k(k-1)/2}$$

$$\le n^{2k} \left( \frac{1}{2} \right) ^{3k(k-1)/2}$$

$$\le n^{(8/3 + 2\epsilon) \log_2{n}} n^{(-3/2)(4/3+\epsilon)((4/3 + \epsilon)\log_2{n}-1)}$$

$$= n^{(8/3 + 2\epsilon) \log_2{n}} n^{-(8/3 + 4 \epsilon + 3\epsilon^2/2) \log_2{n} +2+3\epsilon/2}$$

$$= n^{-(2  + 3\epsilon /2) (\epsilon \log_2{n} - 1)}$$

$$=o(1).$$

\end{proof}

\begin{proof}[Proof of Corollary \ref{corollary 4}] Define the {\it
density} of a graph with $v$ vertices and $e$ edges to be $\lambda
= e/v$. We say a graph is {\it strictly balanced} if the density
of the graph itself is strictly greater than the density of any of
its subgraphs.

Let $H$ be any strictly balanced graph of density $\lambda$. It is
classical that $p=n^{-1 / \lambda}$ is a sharp threshold for
$G(n,p)$ containing $H$ as a subgraph \cite{Bollo}. In particular,
if $p=n^\alpha$ and $\alpha > -1 / \lambda$ then $G(n,p)$ a.a.
contains $H$ as a subgraph, and if $\alpha < -1 / \lambda$ then
$G(n,p)$ a.a. doesn't contain $H$.

Since $K_k$ and $X_k$ are both strictly balanced we may apply this
result twice. The density of $K_k$ is $(k-1)/2$, and the density
of $X_k$ is $3(k-1)/4$. So if $p=n^\alpha$ with $\frac{-2}{k+1}<
\alpha < \frac{-4}{3(k+1)}$, then $G(n,p)$ a.a. contains $K_{k+2}$
but not $X_{k+2}$ subgraphs. This implies that $
\widetilde{H}_k(\| \mathcal{N}[G(n,p)] \|) \neq 0$ by Theorem
\ref{Theorem 3} once we check the detail that at least one of
these $K_{k+2}$ subgraphs can't be extended to a $K_{k+3}$. In
fact a randomly chosen clique will do the job. The conditional
probability that a given $K_{k+2}$ extends to a $K_{k+3}$ is
easily seen to be bounded above by $np^{k+2}$, and $np^{k+2}=o(1)$
since $p=n^\alpha$ with $\alpha<\frac{-4}{3(k+1)}$.

\end{proof}

\section{Connectivity, cliques, and chromatic number}

The chromatic number $\chi(G(n,1/2)$ is tightly concentrated
around $n/\log_2{n}$. For comparison, the clique number is almost
always close to $2\log_2{n}$. As a corollary of what we've shown
here, the connectivity of the neighborhood complex, somewhere
between $\log_2{n}$ and $(4/3)\log_2{n}$, is almost always less
than the clique number.

Similar remarks hold for all monotone functions $p=p(n)$. The
asymptotic picture that emerges is the following. The neighborhood
complex strong deformation retracts to a complex of dimension $d$,
which is $d/4$-connected, with nonvanishing homology between
dimensions $d/3$ and $d/2$, where the clique number is $d/2$. We
see that the connectedness of the neighborhood complex won't do
better than the clique number as a lower bound on chromatic number
for random graphs; the maximal cliques themselves actually
represent nontrivial homology classes.

Recent work of Eric Babson and Dmitry Kozlov
\cite{Babson1,Babson2,Babson3} provides new examples of
topological lower bounds on chromatic number and a more general
setting in which to work. However, Carsten Schultz put bounds on
the strength of these bounds \cite{Schultz}. In particular, the
$\mathbb{Z}_2$-index of the neighborhood complex provides the
strongest known topological bound on chromatic number. This may in
general be higher than the connectivity of the neighborhood
complex. But by what we've shown here, even the
$\mathbb{Z}_2$-index won't do much better for random graphs as a
lower bound on chromatic number than connectivity, since the
dimension of the retract is an upper bound on the index.

\section{Random simplicial complexes and unimodality}

One justification for random graph theory is that it provides
models for ``typical'' graphs. This can be made precise in a few
ways. For example, $G(n,1/2)$ is the uniform distribution on all
graphs on vertex set $[n]$. Any property that $G(n,1/2)$ a.a. has
is a property of almost all graphs. Or for another example, the
Szemer\'edi Regularity Lemma states that every graph is well
approximated by random graphs.

Every neighborhood complex is homotopy equivalent to a free
$\mathbb{Z}_2$-complex, via Lovasz's retract. Up to homotopy, the
converse also holds \cite{Csorba}.

\begin{theorem}[Csorba] Given a finite simplicial complex
$\Delta$ with a free $\mathbb{Z}_2$-action, there exists a graph
$G$ such that $\| \mathcal{N}[G] \|$ is homotopy equivalent to $\|
\Delta \|$.
\end{theorem}

So neighborhood complexes of random graphs asymptotically give a
probability distribution on all finite triangulable
$\mathbb{Z}_2$-spaces as $n \rightarrow \infty$, at least up to
homotopy type.

Little seems to be written so far about random simplicial
complexes. However, Nathan Linial and Roy Meshulam recently
studied $H_1(\| Y(n,p) \|,\mathbb{Z}_2)$ for random
$2$-dimensional simplicial complexes $Y(n,p)$ \cite{Nati}. Their
definition of $Y(n,p)$ is a natural extension of the
Erd\H{o}s-R\'enyi random graph $G(n,p)$; $Y(n,p)$ has vertex set
$[n]$ and edge set ${[n] \choose 2}$, with each $2$-face appearing
independently with probability $p$. One advantage of the
Linial-Meshulam model is that vanishing of homology is a monotone
property. That is, once enough $2$-faces have been added that
$H_1(\| Y(n,p) \|,\mathbb{Z}_2)$ vanishes, adding more $2$-faces
can't ever make it nonvanishing. (This particular fact doesn't
depend on the coefficients of homology, but only on the definition
of random $2$-complex. Simple connectivity is also a monotone
property but it's still not known where the threshold function
lies.)

Most properties of $G(n,p)$ that have been studied to date are
monotone graph properties, in contrast to what we've studied in
this article, where vanishing of homology is clearly not monotone.
But in this setting unimodality seems like a natural substitute
for montonicity.

In another article \cite{Kahle}, we study the {\it clique complex}
of a random graph, which is the simplicial complex with all
complete subgraphs for its faces. The results are analogous to
what we find here, although we also study the expectation of the
Betti numbers. Denote the $k$th Betti number by $\beta_k$. We
conjecture that for both random neighborhood and clique complexes,
for any fixed $k$ and large enough $n$ depending on $k$, the
expectation $E[\beta_k]$ is a unimodal function of $p$.

\section{Acknowledgements}

The author wishes to thank his advisor Eric Babson for
inspiration, guidance, and patience; fellow graduate students,
particularly Anton Dochtermann, Alex Papazoglou, and David Rosoff,
for many helpful conversations; and Sara Billey for generous
support. Any mistakes are his own.

\end{document}